\documentclass{article}
\usepackage{amsmath,amsthm,latexsym,graphicx}

\renewcommand\d{\mathrm{d}}

\newcommand\bc{\mathbf{c}}
\newcommand\bde{\boldsymbol{\delta}}
\newcommand\bE{\mathbf{E}}
\renewcommand\bf{\mathbf{f}}
\newcommand\bF{\mathbf{F}}
\newcommand\bg{\mathbf{g}}
\newcommand\bk{\mathbf{k}}
\newcommand\bL{\mathbf{L}}
\newcommand\bR{\mathbf{R}}
\newcommand\bs{\mathbf{s}}
\newcommand\bT{\mathbf{T}}
\newcommand\by{\mathbf{y}}
\newcommand\byR{\mathbf{y}_{\!\raisebox{-1pt}{$\scriptstyle R$}}}
\newcommand\bz{\mathbf{z}}

\newcommand\hbf{{\hat{\bf}}}

\newcommand\tC{{\tilde{C}}}
\newcommand\ts{{\tilde{s}}}
\newcommand\tit{{\tilde{t}}}

\newcommand\blt{\tau^{\scriptscriptstyle\,|}}

\newcommand\sd{\operatorname{sd}}
\newcommand\sgn{\operatorname{sgn}}

\let\ga\gamma
\let\pa\partial
\let\ss\scriptstyle

\newtheorem{theorem}{Theorem}
\newtheorem{lemma}[theorem]{Lemma}
\newlength\jnheight
\newlength\jndepth
\newsavebox\jnbox
\newcommand\jnstrut[1]{%
  \sbox\jnbox{#1}%
  \settoheight\jnheight{\usebox\jnbox}%
  \addtolength\jnheight\jot%
  \settodepth\jndepth{\usebox\jnbox}%
  \addtolength\jndepth\jot%
  \raisebox{0pt}[\jnheight][\jndepth]{\usebox\jnbox}}
\newcommand\jnstrutb[1]{%
  \sbox\jnbox{#1}%
  \settodepth\jndepth{\usebox\jnbox}%
  \addtolength\jndepth\jot%
  \raisebox{0pt}[\height][\jndepth]{\usebox\jnbox}}

\newcommand\treeA{\begin{picture}(6,5)(-3,-2.5) \put(0,0){\circle*{3}} \end{picture}}
\newcommand\ptreeA{(\treeA)}

\newcommand\treeB{%
  \begin{picture}(4,7)(-2,1) \put(0,0){\circle*{3}} \put(0,7){\circle*{3}} \put(0,0){\line(0,1){7}} \end{picture}}
\newcommand\ptreeB{(\treeB)}

\newcommand\treeC{%
  \raisebox{-4pt}{\begin{picture}(18,14.5)(-9,-2)
    \put(0,0){\circle*{3}} \put(0,0){\line(-2,3){7}} \put(-7,10.5){\circle*{3}} \put(0,0){\line(2,3){7}} 
    \put(7,10.5){\circle*{3}}
  \end{picture}}}
\newcommand\ptreeC{\bigl(\treeC\bigr)}

\newcommand\treeD{%
  \raisebox{-4pt}{\begin{picture}(4,14)(-2,-2)
    \put(0,0){\circle*{3}} \put(0,0){\line(0,1){5}} \put(0,5){\circle*{3}} \put(0,5){\line(0,1){5}}
    \put(0,10){\circle*{3}}
  \end{picture}}}
\newcommand\ptreeD{\bigl(\,\treeD\,\bigr)}

\newcommand\treeE{%
  \raisebox{-4pt}{\begin{picture}(24,14.5)(-12,-2)
    \put(0,0){\circle*{3}} \put(0,0){\line(-1,1){10}} \put(-10,10){\circle*{3}} \put(0,0){\line(0,1){10}}
    \put(0,10){\circle*{3}} \put(0,0){\line(1,1){10}} \put(10,10){\circle*{3}}
  \end{picture}}}
\newcommand\ptreeE{\bigl(\treeE\bigr)}

\newcommand\treeF{%
  \raisebox{-10pt}{\begin{picture}(18,24.5)(-9,-2)
    \put(0,0){\circle*{3}} \put(0,0){\line(-2,3){7}} \put(-7,10.5){\circle*{3}} \put(0,0){\line(2,3){7}}
    \put(7,10.5){\circle*{3}} \put(7,10.5){\line(0,1){10}} \put(7,20.5){\circle*{3}}
  \end{picture}}}
\newcommand\ptreeF{\biggl(\treeF\biggr)}

\newcommand\treeG{%
  \raisebox{-10pt}{\begin{picture}(10,25)(-5,-2)
    \put(0,0){\circle*{3}} \put(0,0){\line(0,1){7}} \put(0,7){\circle*{3}} \put(0,7){\line(0,1){7}}
    \put(0,14){\circle*{3}} \put(0,14){\line(0,1){7}} \put(0,21){\circle*{3}}
  \end{picture}}}
\newcommand\ptreeG{\biggl(\treeG\biggr)}

\newcommand\treeH{%
  \raisebox{-10pt}{\begin{picture}(18,24.5)(-9,-2)
    \put(0,0){\circle*{3}} \put(0,0){\line(0,1){10}} \put(0,10){\circle*{3}} \put(0,10){\line(-2,3){7}}
    \put(-7,20.5){\circle*{3}} \put(0,10){\line(2,3){7}} \put(7,20.5){\circle*{3}}
  \end{picture}}}
\newcommand\ptreeH{\biggl(\treeH\biggr)}

\newcommand\treeZa{%
  \raisebox{-4pt}{\begin{picture}(6,7)(-3,-2) 
    \put(0,0){\circle*{3}} \put(0,8){\makebox[0pt][c]{$\ss\tau$}} \put(0,0){\line(0,1){7}} 
  \end{picture}}}
\newcommand\ptreeZa{\bigl(\treeZa\bigr)}

\newcommand\treeZk{%
  \raisebox{-7pt}{\begin{picture}(26,14.5)(-13,-2)
    \put(0,0){\circle*{3}} \put(0,0){\line(-1,1){10}} \put(-10,12){\makebox[0pt][c]{$\ss\tau_1$}}
    \put(0,12){\makebox[0pt][c]{$\ss\ldots$}} \put(0,0){\line(1,1){10}} \put(10,12){\makebox[0pt][c]{$\ss\tau_k$}}
  \end{picture}}}
\newcommand\ptreeZk{\Bigl(\treeZk\Bigr)}

\begin{document}

\title{A Priori Estimates for the Global Error Committed by Runge--Kutta Methods\\for a Nonlinear Oscillator}
\author{Jitse Niesen\thanks{\ Department of Applied Mathematics and Theoretical Physics, University of Cambridge, Silver Street,
    Cambridge CB3 9EW, United Kingdom, e-mail \texttt{J.Niesen@damtp.cam.ac.uk}.}}
\date{19 December 2001}
\maketitle
\thispagestyle{empty}

\begin{abstract}
  The Alekseev--Gr{\"o}bner lemma is combined with the theory of modified equations to obtain an \emph{a priori} estimate for the
  global error of numerical integrators. This estimate is correct up to a remainder term of order~$h^{2p}$, where $h$ denotes the
  step size and $p$ the order of the method. It is applied to a class of nonautonomous linear oscillatory equations, which
  includes the Airy equation, thereby improving prior work which only gave the $h^p$~term. 
  
  Next, nonlinear oscillators whose behaviour is described by the Emden--Fowler equation \mbox{$y'' + t^\nu y^n = 0$} are
  considered, and global errors committed by Runge--Kutta methods are calculated. Numerical experiments show that the resulting
  estimates are generally accurate. The main conclusion is that we need to do a full calculation to obtain good estimates: the
  behaviour is different from the linear case, it is not sufficient to look only at the leading term, and merely considering the
  local error does not provide an accurate picture either.
\end{abstract}

\section{Introduction}
\label{s:intro}

This paper is about the numerical solution of ordinary differential equations. Unfortunately, the outcome of the calculation is
almost invariably not exact, but carries an error. The error committed by a single step of the numerical method is called the
\emph{local (truncation) error}. But in general, a single step is not sufficient to cover the time interval of interest, so we are
interested in the error at a specific time, after several steps. This quantity is called the \emph{global error}.

It is important to have at least a rough idea of the size of the global error, so we need to \emph{estimate} it. We can
distinguish two approaches. One is to solve the equation numerically, and to use the result somehow to find an \emph{a posteriori}
estimate. The other approach is to use the analytic solution, or at least some knowledge about it. This gives rise to \emph{a
  priori} estimates.

Both types of estimates have their respective strengths. If one has actually computed a numerical solution and wants to know how
far off it is, one probably should use an \emph{a posteriori} estimate. However, for the purpose of comparing different methods,
or that of devising methods which are particularly suited for a certain class of problems, there is often no choice but to use
\emph{a priori} estimates.

Skeel~\cite{skeel:thirteen} wrote a nice overview of various methods to obtain \emph{a posteriori} estimates of the global error,
but here we will concentrate on \emph{a priori} estimates. It is rather hard to find accurate estimates of the global error, as
the following quote by Lambert shows:~\cite[p.~57]{lambert:numerical}
\begin{quote}
  The [local truncation error] and the starting errors accumulate to produce the [global error], but this accumulation process is
  very complicated, and we cannot hope to obtain any usable general expression for the [global error]. However, some insight into
  the accumulation process can be gleaned by looking at an example.
\end{quote}
So the challenge is to find a class of equations which is specific enough to allow the resulting estimates to be of some use, yet
general enough to encompass many problems arising in practice. The class we are aiming for here, can loosely be described as
nonautonomous oscillatory equations.

Among the first \emph{a priori} estimates were those obtained by Henrici~\cite{henrici:discrete} and Gragg~\cite{gragg:repeated},
who expanded the global error in powers of the step size. This expansion was generalized and simplified by Hairer and
Lubich~\cite{hairer.lubich:asymptotic*1}, and used by Cano and Sanz-Serna~\cite{cano.sanz-serna:error*1,cano.sanz-serna:error} to
analyse the behaviour of numerical methods for systems with periodic orbits.

A parallel path of research is based on the Alekseev--Gr\"obner lemma, to be introduced in the next section. We merely point to
the work of Iserles and S\"oderlind~\cite{iserles.soderlind:global}, Calvo and Hairer~\cite{calvo.hairer:accurate}, and
Viswanath~\cite{viswanath:global}. A comparison of both approaches can be found in the paper by Hairer and
Lubich~\cite{hairer.lubich:asymptotic}. 

The recent paper of Iserles~\cite{iserles:on*1} gives an \emph{a priori} estimate for the global error committed for certain linear
oscillators. This estimate only concerns the first term in the expansion of the global error in powers of the step size. The goal
of the research described here, is to extend this work in two directions. Firstly, more (but not all) terms of the expansion are
derived. Numerical experiments will show that it is sometimes necessary to include these terms. Secondly, the effects of
nonlinearity are studied by considering a specific example.

This example is the Emden--Fowler equation $y'' + t^\nu y^n = 0$. For a particular range of parameters, the solutions of this
nonlinear equation are oscillatory. It was selected because it is neither linear nor autonomous, and the oscillatory behaviour of
its solutions will probably lead to substantial cancellations causing more difficulties. Yet it is still amenable to analysis, in
particular because we can derive its asymptotic solution for large~$t$. Finally, it seems a logical sequel to the Airy equation
$y'' + ty = 0$, considered by Iserles~\cite{iserles:on*1}.

The remainder of this paper is organized as follows. The next section introduces the Alekseev--Gr\"obner lemma and the theory of
modified equations, and shows how they can be used to estimate the global error. The resulting \emph{a~priori} estimates are exact
up to a term of order~$h^{2p}$, where $h$ is the step size and $p$ the order of the method. In Section~\ref{s:airy}, we show how
these estimates look like for a class of linear oscillators, which includes the Airy equation. We condense the calculations
because many are already described by Iserles~\cite{iserles:on*1}. The Emden--Fowler equation is introduced in
Section~\ref{s:nonlin}, where we also describe its asymptotic solution. Section~\ref{s:global} contains the actual calculation of
global error estimates. The estimates are supported by numerical experiments, as described in Section~\ref{s:num}. The last
section discusses the results, draws them into comparison with the work done by others, and presents some pointers for further
research.

\section{Theory}
\label{s:theory}

Suppose that we are solving the ordinary differential equation
\begin{equation}
  \label{e:ode}
  \by'(t) = \bf(t,\by(t)), \qquad \by(t_0) = \by_0,
\end{equation}
where $\by(t) \in \bR^d$ is a vector. We assume that $\bf$ is $C^\infty$ to avoid technicalities which are irrelevant to our
argument.

Associated to this differential equation is the \emph{flow map}, which is denoted here by~$\Phi_s^t(\by_0)$, where $s$ stands for
the initial time, $t$~the final time, and $\by_0$~the starting position. The map $\Phi_s^t : \bR^d\to\bR^d$ is defined by
\begin{equation}
  \label{e:flow}
  \Phi_t^t(\by) = \by \quad\text{and}\quad \tfrac{\d}{\d{t}} \Phi_s^t(\by) = \bf \bigl( t, \Phi_s^t(\by) \bigr).
\end{equation}
Its Jacobian matrix will be denoted by $D\Phi_s^t$ and is called the \emph{variational flow}. It plays an important role in this
paper, since it describes how perturbations are propagated by the flow of~\eqref{e:ode}. This follows from the Alekseev--Gr\"obner
lemma (see eg.~\cite{hairer.nrsett.ea:solving}), which extends the variation-of-constants technique to nonlinear equations.
\begin{lemma}[Alekseev--Gr\"obner]
  \label{l:ag}
  Denote by $\by$ and $\bz$ the solutions of
  \begin{subequations}
    \begin{align}
      \by' &= \bf(t,\by),              & \by(t_0) &= \by_0, \label{e:ag_ode1} \\
      \bz' &= \bf(t,\bz) + \bg(t,\bz), & \bz(t_0) &= \by_0, \label{e:ag_ode2}
    \end{align}
  \end{subequations}
  respectively and suppose that $\pa\bf/\pa\by$ exists and is continuous. Then the solutions of~\eqref{e:ag_ode1} and of the
  ``perturbed'' equation~\eqref{e:ag_ode2} are connected by
  \begin{equation}
    \bz(t) = \by(t) + \int_{t_0}^t D\Phi_s^t(\bz(s)) \, \bg\bigl(s,\bz(s)\bigr) \, \d{s}. \label{e:ag}
  \end{equation}
\end{lemma}

Now suppose that one uses a numerical one-step method with fixed step size~$h$ to solve the differential equation~\eqref{e:ode}.
Letting $\by_0=\by(t_0)$, we denote the subsequent results of the method by $\by_1,\by_2,\by_3,\dots$, and set $t_k = t_0 + kh$.

Of course, we hope that numerical results are close to the exact solution, $\by_k \approx \by(t_k)$. The quality of this
approximation is measured by the \emph{global error}, which is defined by $\bE_h(t_k) = \by_k - \by(t_k)$. We say that a method is
of \emph{order}~$p$ if $\bE_h(t) = {\cal O}(h^p)$.

As mentioned in the introduction, the goal of the analysis presented in this paper is to obtain \emph{a priori} estimates for the
global error. The general idea is to view the numerical ``flow'' as a perturbation of the actual flow of the equation. As the
global error is the difference between the numerical and the exact solution, the Alekseev--Gr\"obner lemma can be used to estimate
it.

There is however a problem here: the numerical method provides only values at times~$t_k$, but not a continuous flow. Iserles and
S\"oderlind~\cite{iserles.soderlind:global} interpolated the numerical results to obtain a function defined on the whole time
interval. Instead, we will use the theory of \emph{modified equations}, an idea already present in the work of Calvo and
Hairer~\cite{calvo.hairer:accurate}. We note that the same approach is taken by Hairer and Lubich~\cite{hairer.lubich:asymptotic}.
They proceed to obtain differential equations for the various terms in the expansion of the global error.

The theory of modified equations, analogous to Wilkinson's backward error analysis in Numerical Linear Algebra, tries to find a
differential equation \mbox{$\bz' = \hbf_h(t,\bz)$} whose solution is close to the numerical results. Hairer and
Lubich~\cite{hairer.lubich:life-span} proved that one can arrange for the solution to be exponentially close (in the step
size~$h$) to numerical results. This is more than sufficient to apply the following theorem.

\begin{theorem}
  \label{prop2}
  Let~$\by(t)$ be the solution of the differential equation $\by'=\bf(t,\by)$ and $D\Phi_s^t$ denote its variational flow. Suppose
  that a numerical method of order~$p$ produces the values~$\{\by_k\}$. If the solution~$\bz(t)$ of the modified equation
  \mbox{$\bz' = \hbf_h(t,\bz)$} satisfies $\bz(t_k) - \by_k = {\cal O}(h^{2p})$ for all~$k$, then
  \begin{equation}
    \label{e:prop2}
    \bE_h(t) = \int_{t_0}^t D\Phi_s^t(\by(s)) \, \bde_h\bigl(s,\by(s)\bigr) \, \d{s} + {\cal O}(h^{2p}),
  \end{equation}
  where $\bde_h(t,\by) = \hbf_h(t,\by) - \bf(t,\by)$.
\end{theorem}

\begin{proof}
  If we take the perturbed equation~\eqref{e:ag_ode2} from the Alekseev--Gr\"obner lemma to be the modified equation, we deduce
  $$
  \bz(t) - \by(t) = \int_{t_0}^t D\Phi_s^t(\bz(s)) \, \bde_h\bigl(s,\bz(s)\bigr) \, \d{s}.
  $$
  The left-hand side of this equation is $\bE_h(t) + {\cal O}(h^{2p})$. Furthermore, since the method is of order~$p$, we have
  $\by(t)-\bz(t) = {\cal O}(h^p)$. Finally, it follows from the theory of modified equations that $\bde_h(t,\by) = {\cal O}(h^p)$
  (see eg.~\cite{hairer:numerical}). Together this proves~\eqref{e:prop2}.
\end{proof}

Note that the constant implied by the ${\cal O}(h^{2p})$~term in~\eqref{e:prop2} depends on~$t$. Hence the above result is only
valid on bounded time intervals. The same goes for the results we will obtain later (Theorems~\ref{prop3}, \ref{th:linosc}
and~\ref{th:ef}). 

\enlargethispage{2.3ex} 
To compute the error estimate~\eqref{e:prop2}, we need to know the exact solution~$\by(t)$, the variational flow~$D\Phi_s^t$, and
the modified equation to compute~$\bde_h$. Note that these can all be determined \emph{a priori}, hence the above theorem
gives an \emph{a priori} estimate for the global error.

In fact, if one is interested in the long-term behaviour of the global error, it is not essential to know the exact solution over
the whole time interval. It often suffices to know only the asymptotic behaviour of the solution for large~$t$, because this will
allow one to compute the dominant contribution to the integral~\eqref{e:prop2}. This approach is actually taken in the
computations that will follow in the next sections.

From now on, we restrict our attention to \emph{Runge--Kutta methods}. We assume that the reader is familiar with the theory of
Runge--Kutta methods and B-series, as explained in the text books~\cite{hairer.nrsett.ea:solving} and~\cite{lambert:numerical},
amongst others. A B-series is a formal power series of the form
\begin{equation}
  \label{e:Bseries}
  B(a,\by) = a(\emptyset) \by + \sum_{\tau\in\bT} \frac{h^{\rho(\tau)}}{\rho(\tau)!} \, \alpha(\tau) \, a(\tau) \, \bF(\tau)(\by).
\end{equation}
Here $\bT$~denotes the set of rooted trees, $a:\bT\to\bR$ is a coefficient function, $\rho(\tau)$~is the order of the tree~$\tau$,
$\alpha(\tau)$ denotes the number of monotone labellings of~$\tau$, and $\bF(\tau)(\by)$ is the elementary differential of the
function~$\bf$ associated with~$\tau$. The result of any Runge--Kutta method can be written as a B-series, and the coefficient
function~$a$ depends on the coefficients of the method.

Hairer~\cite{hairer:numerical} uses the concept of B-series to derive a formula for the modified equation, taking an approach
similar to that of Benettin and Giorgilli~\cite{benettin.giorgilli:on}. Suppose that the outcome of the Runge--Kutta method is
described by the B-series with coefficient function~$a$. Then the modified equation is given by $\bz' = \frac1h B(b,\bz)$, where
the coefficients are recursively defined by
\begin{equation}
  \label{e:b-def}
  b(\emptyset)=0,\, b\ptreeA=1, \text{ and }
  b(\tau) = a(\tau) - \sum_{j=2}^{\rho(\tau)} \frac{1}{j!} \pa_b^{j-1} b(\tau) \text{ for } \rho(\tau)\ge2.
\end{equation}
In this formula, $\pa_b$~is a Lie derivative, defined as follows.  Given two coefficient functions~$b$ and~$c$
with~$b(\emptyset)=0$, the Lie derivative of~$B(c,\by)$ with respect to the vector field~$B(b,\by)$ is defined as $\frac\d{\d{t}}
B(c,\by_b(t))$, where $\by_b(t)$ is a (formal) solution of the differential equation $\by'(t) = B(b,\by(t))$. It turns out that
this is again a B-series, and its coefficient function is denoted by~$\pa_b c$. An explicit formula for~$\pa_b c$ is given
in~\cite{hairer:numerical}.

We now use this expression for the modified equation to rewrite the \emph{a priori} estimate given in Theorem~\ref{prop2}.

\begin{theorem}
  \label{prop3}
  Let~$\by(t)$ be the solution of the differential equation $\by'=\bf(t,\by)$ and let $D\Phi_s^t$ denote its variational flow.
  Suppose that this equation is solved by a Runge--Kutta method with B-series coefficient function~$a$. Let~$b$ be defined as
  in~\eqref{e:b-def}. If the numerical method has order~$p$, then the global error satisfies:
  \begin{equation}
    \label{e:prop3}
    \begin{aligned}
      {}& \bE_h(t) = \sum_{\tau\in\bT_{\geq2}} h^{\rho(\tau)-1} b(\tau) \frac{\alpha(\tau)}{\rho(\tau)!} {\cal I}_\tau(t) 
      + {\cal O}(h^{2p}), \\[-1.1ex] 
      {}& \text{where } {\cal I}_\tau(t) = \int_{t_0}^t D\Phi_s^t(\by(s)) \bF(\tau)(\by(s)) \,\d s. 
    \end{aligned}
  \end{equation}
  Here $\bT_{\geq2}$ denotes the set of trees with order at least~2.
\end{theorem}

\begin{proof}
  As mentioned above, the right-hand side of the modified equation is $\frac1h B(b,\bz)$, with~$b$ as in~\eqref{e:b-def}. Now
  $\bF\ptreeA=\bf$, so the first term in this B-series is the right-hand side of the original equation~\eqref{e:ode}. Hence their
  difference $\bde_h(t,\by)$ equals
  $$
  \sum_{\tau\in\bT_{\ge2}} \frac{h^{\rho(\tau)-1}}{\rho(\tau)!} \, \alpha(\tau) \, b(\tau) \, \bF(\tau)(\by).
  $$
  We now substitute this expression in~\eqref{e:prop2} and move the scalar factors out of the integral (remember that the
  variational flow~$D\Phi_s^t$ is a linear operator). This results in~\eqref{e:prop3}.
\end{proof}

Note that the error estimate~\eqref{e:prop3} nicely separates the problem and the method. The method only enters the estimate via
the coefficients~$b(\tau)$. On the other hand, the value of the integrals ${\cal I}_\tau$ is completely determined by the
particular equation one is solving. As their role in the global error estimate~\eqref{e:prop3} is similar to the role of the
elementary differential~$\bF(\tau)(\by)$ in the local error, we will call ${\cal I}_\tau$ the \emph{elementary integral}
associated with~$\tau$.

\section{The Airy equation and other linear oscillators}
\label{s:airy}

Iserles~\cite{iserles:on*1} studied in detail the leading term of the global error committed by various methods, including
Runge--Kutta methods, for the Airy equation
\begin{equation}
  \label{e:airy}
  y'' + ty = 0, \quad t\geq0, \qquad y(0)=y_0, \quad y'(0)=y'_0.
\end{equation}
Here we will use Theorem~\ref{prop3} to extend his analysis to include higher-order terms. This will prepare us for the study of
nonlinear variants of the Airy equation, viz.~\eqref{e:ef}.

The asymptotic solution of~\eqref{e:airy} for large~$t$ can be found by the Liouville--Green
approximation~\cite{olver:asymptotics}, which is also known, particularly to physicists, as the WKB(J)-approximation:
\begin{equation}
  \label{e:a-assol}
  \by(t) \approx \Lambda(t) R(\theta(t)) \bs_0\text{, where } \by(t) = \begin{bmatrix} y(t) \\ y'(t) \end{bmatrix}.
\end{equation}
Here $\bs_0$ is a vector whose value is determined by the initial conditions, and $\theta(t) = \tfrac23 t^{3/2}$. Furthermore,
$$
\Lambda(t) = \begin{bmatrix} t^{-1/4} & 0 \\ 0 & t^{1/4} \end{bmatrix} \quad\text{and}\quad
R(\theta) = \begin{bmatrix} \cos\theta & \sin\theta \\ -\sin\theta & \cos\theta \end{bmatrix}.
$$
Note that $\Lambda(t)$ is a scaling matrix and $R(\theta)$ is a rotation.

We want to apply~\eqref{e:prop3} to derive an estimate for the global error, so we have to calculate the elementary
differentials~$\bF(\tau)(\by)$. This has already been done by Iserles~\cite{iserles:on*1}. He concludes that amongst all the trees
of order~$\rho$, the tall tree without branches dominates. We denote this tree by~$\blt_\rho$ (anticipating our analysis of the a
nonlinear variant of~\eqref{e:airy}, compare with the trees in Table~\ref{t:trees}). The elementary differential associated to
this tree is $\bF(\blt_\rho)(\by) = \hat{F}_{\!\rho}\by$, with
\begin{equation}\label{e:a-hatf} 
\hat{F}_{2r} \approx \begin{bmatrix} (-t)^r & 0 \\ 0 & (-t)^r \end{bmatrix} \quad\text{and}\quad
\hat{F}_{2r+1} \approx \begin{bmatrix} 0 & (-t)^r \\ (-t)^{r+1} & 0 \end{bmatrix}.
\end{equation}
The next step is to compute the elementary integral~${\cal I}_{\blt_\rho}$, defined in~\eqref{e:prop3}. We will use the
shorthand~${\cal I}_\rho$ for this elementary integral. Hence we need the variational flow~$D\Phi_s^t$. But since the Airy
equation~\eqref{e:airy} is linear, the variational flow equals the flow map~$\Phi_s^t$, so we can deduce from~\eqref{e:a-assol}
that
\begin{equation}\label{e:a-dphi}
D\Phi_s^t = \Phi_s^t \approx \Lambda(t) R \bigl( \theta(t)-\theta(s) \bigr) \Lambda^{-1}(s).
\end{equation}
Now we need to substitute~\eqref{e:a-assol}, \eqref{e:a-hatf}, and~\eqref{e:a-dphi} in the elementary integral \mbox{${\cal
    I}_\rho = \int_0^t D\Phi_s^t \hat{F}_{\!\rho} \by(s) \,\d{s}$}. For even~$\rho$ we find
\begin{align*}
  {\cal I}_{2r} &\approx \int_0^t \Lambda(t) R \bigl( \theta(t)-\theta(s) \bigr) \Lambda^{-1}(s) 
  \cdot (-s)^r \Lambda(s) R(\theta(s)) \bs_0 \,\d s \\
  &= \int_0^t (-s)^r \,\d s \cdot \Lambda(t)R(\theta(t))\bs_0 = (-1)^r \frac{t^{r+1}}{r+1} \by(t). 
\end{align*}
For odd~$\rho$, the calculation is only a bit more complicated. First we compute the matrix product that appears in the middle of
the integrand:
$$
\Lambda^{-1}(s) \hat{F}_{2r+1} \Lambda(s) = (-1)^r s^{r+\frac12} R(\tfrac12\pi).
$$
Note that $R(\tfrac12\pi) = \bigl[\begin{smallmatrix}0&1\\-1&0\end{smallmatrix}\bigr]$. With this result, we can tackle the
integral,
$$  
{\cal I}_{2r+1} \approx \!\int_0^t (-1)^r s^{r+\frac12} \,\d{s} \cdot \Lambda(t) R \bigl( \theta(t)+\tfrac12\pi \bigr) \bs_0 
= (-1)^r \frac{t^{r+\frac32}}{r+\tfrac32} \byR(t).
$$
Here $\byR(t)$ is the (asymptotic) solution of the Airy equation~\eqref{e:airy} with opposite phase. It is given by
\begin{equation}
  \label{e:yR}
  \byR(t) = \Lambda(t) R \bigl( \theta(t)+\tfrac12\pi \bigr) \bs_0.
\end{equation}
To find an estimate for the error, we have to substitute the above integrals in~\eqref{e:prop3}. This gives
\begin{multline*}
  \bE_h(t) \approx C_1ht^2\by(t) + C_2h^2t^{5/2}\byR(t) + C_3h^3t^3\by(t) \\[-1ex]
  + C_4h^4t^{7/2}\byR(t) + \cdots + {\cal O}(h^{2p}), \text{where } C_k = \frac{b(\blt_{k+1})}{(k+1)!\,(\frac12k+1)}.
\end{multline*}
Note that for a method of order~$p$, the B-series coefficient of any tree whose order does not exceed~$p$ vanishes. Therefore,
the $h^k$~term vanishes as well for $k<p$, in agreement with the general theory. 

Note furthermore that the nature of the particular Runge--Kutta method used, has no influence on the growth rate of the global
error. It does not matter whether the method is for instance A-stable or symplectic. There is one exception to this, as pointed
out by Orel~\cite{orel:runge-kutta}: for methods with $b(\blt_{p+1}) = 0$, the dominant term in the above estimate vanishes. Hence
such methods should perform better than general Runge--Kutta methods.

It is possible to generalize the asymptotic solution~\eqref{e:a-assol} to a class of nonautonomous oscillators. Under the
assumptions given by Iserles~\cite{iserles:on*1} and reproduced in Theorem~\ref{th:linosc} below, the asymptotic solution of the
equation $y'' + g(t)y = 0$ is still given by~\eqref{e:a-assol}, except that now $\theta(t) = \int_0^t \sqrt{g(s)} \,\d{s}$. The
rest of the computation stays essentially the same. Its result is given by the following theorem.

\begin{theorem}
  \label{th:linosc}
  Suppose we are solving the nonautonomous linear oscillator
  $$
  y'' + g(t)y = 0, \quad t\geq0, \qquad y(0)=y_0, \quad y'(0)=y'_0,
  $$
  where $g(t)$ is positive for sufficiently large~$t$, $|g^{(\ell)}(t)| = o\bigl(g(t)^{1/\ell}\bigr)$ as $t\to\infty$, and
  $\limsup_{t\to\infty} g(t) = +\infty$. If we use a Runge--Kutta method with B-series coefficient function~$a$, and $b$~is
  defined as in~\eqref{e:b-def}, then the global error is
  \begin{multline*}
    \bE_h(t) \approx \left[ \sum_{p\le2r<2p} (-1)^r \frac{b(\blt_{2r+1})}{(2r+1)!} \,h^{2r} 
      \int_{t_0}^t g(s)^{r+\frac12} \,\d{s} \right] \byR(t) \\
    + \left[ \sum_{p\le2r+1<2p} (-1)^{r+1} \frac{b(\blt_{2r+2})}{(2r+2)!} \,h^{2r+1} 
      \int_{t_0}^t g(s)^{r+1} \,\d{s} \right] \by(t) + {\cal O}(h^{2p}).
  \end{multline*}
  Here $p$ denotes the order of the method, and $\by(t)$ and $\byR(t)$ are given by~\eqref{e:a-assol} and~\eqref{e:yR}
  respectively, with $\theta(t) = \int_0^t \sqrt{g(s)} \,\d{s}$. \qed%
\end{theorem}

\section{The Emden--Fowler equation}
\label{s:nonlin}

The previous section should be considered as a mere warm-up to the real challenge: nonlinear oscillators. We will consider
nonautonomous, nonlinear oscillators whose behaviour is described by the Emden--Fowler equation
\begin{equation}
  \label{e:ef}
  y'' + t^\nu y^n = 0, \quad t\geq0, \qquad y(0)=y_0, \quad y'(0)=y'_0.
\end{equation}
This equation with $\nu=1-n$ (in which case it is commonly called the Lane--Emden equation) was originally proposed to model stars.
Considering the star as a radially symmetric gaseous polytrope of index~$n$ in thermodynamic and hydrostatic equilibrium, the
relation between the density and the distance to the centre satisfies the Lane--Emden equation~\cite{chandrasekhar:introduction}.
The Emden--Fowler equation also arises in the fields of gas dynamics, fluid mechanics, relativistic mechanics, nuclear physics,
and the study of chemically reacting systems (see~\cite{wong:on} and references therein).

Note that the choice $\nu=n=1$ reduces the Emden--Fowler equation to the Airy equation studied in the previous section. From now
on we will assume that $n$ is an odd integer above~$1$, and that $\nu > -\tfrac12(n+3)$. These conditions assure that oscillatory
solutions exists, as explained in the survey by Wong~\cite{wong:on}. We remark incidentally that this remains true for any real
$n>1$, if we replace the equation~\eqref{e:ef} by $y'' + t^\nu \sgn(y) \, |y|^n = 0$, where $\sgn(y)$ denotes the sign of~$y$.
However, we then lose analyticity at $y=0$.

The equation~\eqref{e:ef} has an obvious scaling symmetry, which can be used to reduce the order. Here we will use a different
approach though, because the asymptotic solution as $t\to\infty$ will suffice for our purposes. From now on, we set
$$
\gamma = \frac{\nu}{n+3}.
$$
The conditions on~$n$ and~$\nu$ imply that $\gamma > -\tfrac12$. Now consider the transformation given by
\begin{equation}
  \label{e:ef-tf}
  y(t) = (1+2\gamma)^{2/(n-1)} t^{-\gamma} u(t^{1+2\gamma}).
\end{equation}
If we apply~\eqref{e:ef-tf} to the differential equation~\eqref{e:ef}, we get
\begin{multline*}
  (1+2\ga)^{2n/(n-1)} t^{3\ga} \bigl( u''(t^{1+2\ga}) + u^n(t^{1+2\ga}) \bigr) \\
  + \ga(1+\ga)(1+2\ga)^{2/(n-1)} t^{-\ga-2} u(t^{1+2\ga}) = 0.
\end{multline*}
For large~$t$, we can (hopefully) neglect the last term as $\ga > -\tfrac12$, and the above equation reduces to $u''+u^n=0$. The
expression $\frac1{n+1}u^{n+1} + \frac12(u')^2$ is an invariant of this equation; it is only an adiabatic invariant to the
original equation. We will denote the solution of $u''+u^n=0$ which satisfies the initial conditions $u(0)=0$ and $u'(0)=1$
by~$w_n(t)$, and note for further reference that this is an odd and periodic function. The general solution of $u''+u^n=0$ is then
given by $u(t) = c_1^{2/(n+1)} w_n(c_1 t+c_2)$. Note that $c_1$ determines the amplitude of the oscillation, while $c_2$
determines the phase. In other, more sophisticated terms, $(c_1,c_2)$ are action-angle coordinates of the Hamiltonian system
corresponding to the Emden--Fowler oscillator (see eg.~\cite{abraham.marsden:foundations}).

It follows that the solution of the Emden--Fowler equation~\eqref{e:ef} is asymptotically given by
\begin{equation}
  \label{e:ef-sol}
  y(t) \approx (1+2\ga)^{2/(n-1)} c_1^{2/(n-1)} t^{-\ga} w_n(c_1t^{1+2\ga}+c_2).
\end{equation}
We repeat that we assumed that $n>1$ is an odd integer and that $\nu > -\tfrac12(n+3)$.

\section{Global error of Runge--Kutta methods}
\label{s:global}

In this section, we repeat the computation of Section~\ref{s:airy} for the nonlinear oscillator~\eqref{e:ef}. We will assume that
the asymptotic solution~\eqref{e:ef-sol} is in fact exact. The numerical experiments reported in Section~\ref{s:num} will show
that this is a valid approximation.

The first task is the computation of elementary differentials. For this we need to convert the differential equation to a system
of autonomous first-order equations,
\begin{equation}
  \label{e:c-fo}
  \begin{array}{l}
      y'_1 = y_2, \\
      y'_2 = -y_3^\nu y_1^n, \\
      y'_3 = 1.
  \end{array}
\end{equation}
Here $y_1$, $y_2$, and $y_3$ correspond to $y$, $y'$, and $t$ in the original equation~\eqref{e:ef} respectively.

It follows that the first component of the elementary differential~$\bF(\tau)(\by)$ satisfies the following recurrence relations
(where the argument~$\by$ is deleted)
$$
F_1\ptreeA = y_2, \qquad F_1\ptreeZa = F_2(\tau), \qquad F_1\ptreeZk = 0\quad(k\ge2).
$$
The second component satisfies
\begin{align*}
  F_2\ptreeA &= -y_3^\nu y_1^n, \\
  F_2\ptreeZk &= -\frac{n!}{(n-k)!} y_3^\nu y_1^{n-k}F_1(\tau_1)\ldots F_1(\tau_k) && (k\le n), \\
  F_2\ptreeZk &= 0 && (k\ge n+1),
\end{align*}
where the derivatives with respect to~$y_3$ are neglected, because they lead to terms which are smaller by a factor~$t^{1+2\ga}$.
Finally, the third component~$F_3$ always vanishes, except for~$F_3\ptreeA$ which equals one. Hence we will drop this component
from now on.

It follows from the recursion relations that the first component of the elementary differential associated with the tree~$\tau$
can be computed as follows. Each vertex contributes a factor to~$F_1(\tau)$. Define the height of a vertex to be the distance to
the root. Then for vertices whose height is even, the factor is $y_2$ if the vertex is a leaf, $1$~if it has one child and zero
otherwise, killing the whole elementary differential.  For vertices with odd height, a vertex with $d$ children contributes a
factor $t^\nu y_1^{n-d}$ (apart from a constant) if $d\le n$, and zero otherwise. So $F_1(\tau)$ vanishes if any vertex with even
height has more than one child, or any vertex with odd height has more than $n$~children. If this is not the case, let $d'$ denote
the number of vertices with odd height, and $\rho$ the total number of vertices, ie.~the order of the tree.  Then there are
$\rho-d'$~vertices with even height, who have $d'$~children, namely the $d'$~vertices with odd height. Since every vertex with
even height has at most one child, we have $\rho-2d'$~leaves contributing a factor~$y_2$ each and $d'$~vertices with exactly one
child, contributing a factor~$1$.  So the total contribution of all the $d'$~vertices with even height is $y_2^{\rho-2d'}$.
Similarly, the $d'$~vertices with odd height have $\rho-d'-1$~children (the root is not a child of any vertex).  So they
contribute a factor $t^{d'\nu} y_1^{nd' - (\rho-d'-1)}$. We conclude that
$$
F_1(\tau)(\by(t)) = C_{1,\tau} t^{d'\nu} y_1^{(n+1)d' - \rho + 1}(t) \, y_2^{\rho-2d'}(t).
$$
By a similar reasoning, we find that
$$
F_2(\tau)(\by(t)) = C_{2,\tau} t^{(\rho-d')\nu} y_1^{n\rho-(n+1)d'}(t) \, y_2^{2d'-\rho+1}(t).
$$
Note that the conditions on the tree which make sure that the elementary differential does not vanish, imply that all the
exponents in the above formulae are positive. Also, the only trees for which both~$F_1$ and~$F_2$ are nonzero, are the tall
branchless trees~$\blt_\rho$.

If we substitute the approximate solution~\eqref{e:ef-sol}, we find that the elementary differentials are
\begin{equation}
  \label{e:ef-eldif}
  \bF(\tau)(\by(t)) \approx \begin{bmatrix}
    C_{3,\tau} t^{\ga(2\rho-1)} w_n^{(n+1)d'-\rho+1}(\tit) \, {w'_n}^{\rho-2d'}(\tit) \\
    C_{4,\tau} t^{\ga(2\rho+1)} w_n^{n\rho-(n+1)d'}(\tit) \, {w'_n}^{2d'-\rho+1}(\tit)
  \end{bmatrix}
\end{equation}
where $\tit = c_1t^{1+2\ga}+c_2$. The growth rate of the elementary differential is determined by the exponent of~$t$. Note that
the variable~$d'$ does not enter in this exponent. The surprising conclusion is that all trees with the same order contribute a
term with the same growth rate, independent of their shape. This is in stark contrast to the linear Airy equation, where the
differential corresponding to the tall tree~$\blt_\rho$ dominates.

The next step is to calculate the elementary integral~${\cal I}_\tau$. For this, we need to multiply the above differential with
the variational flow matrix and integrate the resulting expression. To compute the variational flow, we introduce the
map~$X_t:\bR^2\to\bR^2$, defined by
\begin{equation}
  \label{e:Xt}
  X_t(\bc) = \begin{bmatrix}
    (1+2\ga)^{2/(n-1)}   c_1^{2/(n-1)}   t^{-\ga} w_n(c_1t^{1+2\ga}+c_2) \\
    (1+2\ga)^{1+2/(n-1)} c_1^{1+2/(n-1)} t^{\ga}  w'_n(c_1t^{1+2\ga}+c_2) 
  \end{bmatrix}.  
\end{equation}
So $X_t$ maps the parameter space to the solution space at time~$t$. Neglecting lower-order terms, it follows that the flow map
satisfies $\Phi_s^t = X_t \circ X_s^{-1}$. Hence we can write the elementary integral~\eqref{e:prop3} as:
\begin{equation}
  \label{e:ef-elint-d}
  {\cal I}_\tau(t) \approx DX_t(\bc) \int_0^t DX_s^{-1}(\by(s)) \, \bF(\tau)(\by(s)) \,\d s. 
\end{equation}
To find the integrand in the above expression, we multiply the inverse of the Jacobian matrix of~\eqref{e:Xt}
with~\eqref{e:ef-eldif}. The result is
\begin{equation}
  \label{e:ig}
  \begin{bmatrix}
    s^{2\ga\rho} \Bigl( C_{5,\tau} w_n^{(n+1)(d'+1)-\rho} (w'_n)^{\rho-2d'} 
    + C_{6,\tau} w_n^{n\rho-(n+1)d'} (w'_n)^{2d'-\rho+2} \Bigr) \\
    s^{2\ga\rho+2\ga+1} \Bigl( C_{7,\tau} w_n^{(n+1)(d'+1)-\rho} (w'_n)^{\rho-2d'} 
    + C_{8,\tau} w_n^{n\rho-(n+1)d'} (w'_n)^{2d'-\rho+2} \Bigr)
  \end{bmatrix},
\end{equation}
where the functions~$w_n$ and~$w'_n$ are evaluated at $\ts = c_1s^{1+2\ga}+c_2$. 

In the calculation, we used that 
$$
w''_n(t) = -w_n^n(t) \quad\text{and}\quad {w'_n}^2(t) = -\tfrac{1}{n+1}w_n^{n+1}(t) + 1.
$$
The first equality is indeed the definition of~$w_n$, and the second one follows by multiplying the first one by $w'_n(t)$ and
integrating, using the initial conditions~$w_n(0)=0$ and $w'_n(0)=1$.

The next step is to integrate~\eqref{e:ig}. But consider the exponents of~$w$ and~$w'$. If $\rho$ is even, these exponents are
also even and hence the integrand is nonnegative. Now recall that $w$ is odd and periodic. So in the case when $\rho$ is odd, the
functions~$w$ and~$w'$ are raised to an odd power, which means that the integrand oscillates around zero. Thus we can expect
cancellations in the latter case, but not if $\rho$ is even. We stress that this phenomenon does not occur in the linear case,
analysed in Section~\ref{s:airy}.

In fact, we have
$$
\int_0^t w_n^\ell(s) \, {w'_n}^m(s) \,\d{s} = \begin{cases}
  \tC_{\ell mn}(t), & \text{if either $\ell$ or $m$ is odd}, \\
  C_{\ell mn}t + \tC_{\ell mn}(t), & \text{if both $\ell$ and $m$ are even}.
\end{cases}
$$
Here $\tC_{\ell mn}(t)$ denotes a oscillatory function with the same period as~$w_n(t)$, and $C_{\ell mn}$ is a constant.  It
follows that
$$
\int_0^t s^k w_n^\ell(\ts) \, {w'_n}^m(\ts) \,\d{s} = 
\begin{cases}
  \tC_{k\ell mn}(\tit) \, t^{k-2\ga} + {\cal O}(t^{k-4\ga-1}), & \text{if $\ell$ or $m$ are odd}, \\
  C_{k\ell mn}t^{k+1} + {\cal O}(t^{k-2\ga}), & \text{if $\ell$ and $m$ are even},
\end{cases}
$$
where again $\tC_{k\ell mn}$ and $C_{k\ell mn}$ denote a periodic function and a constant, respectively.

We can use this result to integrate~\eqref{e:ig}, and find
\begin{equation}
  \label{e:ef-int}
  \int_0^t DX_s^{-1} \bF(\tau)(\by) \,\d{s} =
  \begin{cases}
    \,\begin{bmatrix} 
      \tC_{9,\tau}(\tit) \, t^{2\ga\rho-2\ga} + {\cal O}(t^{2\ga\rho-4\ga-1}) \\
      \tC_{10,\tau}(\tit) \, t^{2\ga\rho+1} + {\cal O}(t^{2\ga\rho-2\ga})
    \end{bmatrix}, & \text{if $\rho$ odd}, \\[3ex]
    \,\begin{bmatrix} 
      C_{9,\tau} t^{2\ga\rho+1} + {\cal O}(t^{2\ga\rho-2\ga}) \\
      C_{10,\tau} t^{2\ga\rho+2\ga+2} + {\cal O}(t^{2\ga\rho+1})
    \end{bmatrix}, & \text{if $\rho$ even}.
  \end{cases}
\end{equation}
To compute the elementary integral~${\cal I}_\tau$, we need to premultiply the integral~\eqref{e:ef-int} with~$DX_t$,
cf.~\eqref{e:ef-elint-d}. But the expression~\eqref{e:ef-int} has an interpretation by itself. Remember that $X_t$ maps the
parameter space to the solution space. So the integral~\eqref{e:ef-int} represents the error in the parameter space. We conclude
that the energy error associated with the tree~$\tau$ grows as~$t^{2\ga\rho+1}$ if $\rho$ is even, and as~$t^{2\ga\rho-2\ga}$ if
$\rho$ is odd. The second component gives the phase error.

Multiplying the Jacobian matrix of the map~$X_t$ with the integral~\eqref{e:ef-int} gives us the elementary integrals
\begin{equation*}
  \label{e:c-elint}
  {\cal I}_\tau(t) = 
  \begin{cases}
    \,\begin{bmatrix} 
      \tC_{11,\tau}(\tit) \, t^{2\ga\rho-\ga+1} + {\cal O}(t^{2\ga\rho-3\ga}) \\
      \tC_{12,\tau}(\tit) \, t^{2\ga\rho+\ga+1} + {\cal O}(t^{2\ga\rho-\ga})
    \end{bmatrix}, & \text{if $\rho$ odd}, \\[3ex]
    \,\begin{bmatrix} 
      \tC_{11,\tau}(\tit) \, t^{2\ga\rho+\ga+2} + {\cal O}(t^{2\ga\rho-\ga+1}) \\
      \tC_{12,\tau}(\tit) \, t^{2\ga\rho+3\ga+2} + {\cal O}(t^{2\ga\rho+\ga+1})
    \end{bmatrix}, & \text{if $\rho$ even}.
  \end{cases}
\end{equation*}
Finally we can find an estimate for the global error by adding the contributions of all trees according to~\eqref{e:prop3}. For
the first component, we find 
$$
E_h(t) \approx \tC_1(\tit)ht^{5\ga+2} + \tC_2(\tit)h^2t^{5\ga+1} + \tC_3(\tit)h^3t^{9\ga+2} + \tC_4(\tit)h^4t^{9\ga+1} + \cdots +
{\cal O}(h^{2p}).
$$
More formally, we have the following result.
\begin{theorem}
  \label{th:ef}
  Suppose one is solving the Emden--Fowler equation $y''+t^\nu y^n=0$, with $\nu>-\frac12(n+3)$ and $n$ an odd integer greater
  than~$1$, with a numerical method of order~$p$ that can be expressed as a B-series. Then the global error has the form
  \begin{multline*}
    \bE_h(t) = \sum_{p\le2r<2p} h^{2r} \begin{bmatrix} 
      \tC_{2r}^1(\tit)\,t^{4\ga r+\ga+1} + {\cal O}(t^{4\ga r-\ga}) \\ 
      \tC_{2r}^2(\tit)\,t^{4\ga r+3\ga+1} + {\cal O}(t^{4\ga r+\ga}) 
    \end{bmatrix} \\
    + \sum_{p\le2r+1<2p} h^{2r+1} \begin{bmatrix} 
      \tC_{2r+1}^1(\tit)\,t^{4\ga r+5\ga+2} + {\cal O}(t^{4\ga r+3\ga+1}) \\ 
      \tC_{2r+1}^2(\tit)\,t^{4\ga r+7\ga+2} + {\cal O}(t^{4\ga r+5\ga+1}) 
    \end{bmatrix} 
    + {\cal O}(h^{2p}).
  \end{multline*}
  Here $\tC^i_k(\tit)$~denotes a function periodic in $\tit = c_1t^{4/3}+c_2$, and $\ga=\nu/(n+3)$. \qed
\end{theorem}

We would again like to draw the reader's attention to the difference between the even and the odd powers of~$h$. Unfortunately, we
do not know what causes this phenomenon, and whether it also occurs for other differential equations.

\section{Numerical experiments}
\label{s:num}

The purpose of this section is to supplement the calculations of the previous section with some numerical experiments. In
particular, we want to see whether we were justified in using the asymptotic solution~\eqref{e:ef-sol}. Furthermore, it could be
interesting to study whether the ${\cal O}(h^p)$~term of the global error~$\bE_h(t)$ dominates, or whether other terms have to be
taken into account.

All experiments are performed with the parameters $n=3$ and $\nu=1$, so the differential equation we are solving is 
\begin{equation}
  \label{e:cairy}
  y''+ty^3=0.
\end{equation}
The reason for this particular choice is that the function $w_n(t)$, the solution of $u''+u^n=0$ satisfying $u(0)=0$ and $u'(0)=1$
(cf.~Section~\ref{s:nonlin}), can be expressed in terms of Jacobi elliptic functions (see eg.~\cite{neville:jacobian}). In fact,
we have $w_3(t) = \sd(t\,|\,\frac12)$. The parameter~$\frac12$ will be dropped from now on. As a consequence, we can calculate the
elementary integral associated with any given tree explicitly. For the first couple of trees, this yields the results listed in
Table~\ref{t:trees}.

\begin{table}
  \begin{center}\begin{tabular}{|l@{}c|r@{${}={}$}lc|} 
    \hline
    \multicolumn{2}{|c|}{Tree} & \multicolumn{2}{c}{Elementary differential} & Elementary integral \\
    \multicolumn{2}{|c|}{$\tau$} & \multicolumn{2}{c}{$\bF(\tau)(\by)$} 
    & \jnstrutb{${\cal I}_\tau(t) = \int_0^t D\Phi_s^t \bF(\tau)(\by) \,\d{s}$}
    \\\hline\hline
    $\blt_2$: & \treeB 
    & $\begin{bmatrix} -y_1^3y_3 \\ -3y_1^2y_2y_3-y_1^3 \end{bmatrix}$
    & $\begin{bmatrix} {\cal O}(t^{1/2}) \\ {\cal O}(t^{5/6}) \end{bmatrix}$ 
    & \jnstrut{$\begin{bmatrix} 
        -\frac{128}{675}\sqrt{2}\,c_1^4\chi t^{17/6}\sd'(\tit) \\[\jot] 
        \frac{256}{2025}\sqrt{2}\,c_1^5\chi t^{19/6}\sd^3(\tit) 
      \end{bmatrix}$}
    \\\hline\hline
    $\tau_3^a$: & \treeC 
    & $\begin{bmatrix} 0 \\ -6y_1y_2^2y_3-6y_1^2y_2 \end{bmatrix}$
    & $\begin{bmatrix} 0 \\ {\cal O}(t^{7/6}) \end{bmatrix}$ 
    & \jnstrut{$\begin{bmatrix} 
        -\frac{32}{135}\sqrt{2}\,c_1^4\chi t^{11/6}\sd'(\tit) \\[\jot] 
        \frac{64}{405}\sqrt{2}\,c_1^5\chi t^{13/6}\sd^3(\tit) 
      \end{bmatrix}$}
    \\\hline
    $\blt_3$: & \treeD 
    & $\begin{bmatrix} -3y_1^2y_2y_3-y_1^3 \\ 3y_1^5y_3^2 \end{bmatrix}$
    & $\begin{bmatrix} {\cal O}(t^{5/6}) \\ {\cal O}(t^{7/6}) \end{bmatrix}$ 
    & \jnstrut{$\begin{bmatrix} 
        -\frac{208}{135}\sqrt{2}\,c_1^4\chi t^{11/6}\sd'(\tit) \\[\jot] 
        \frac{416}{405}\sqrt{2}\,c_1^5\chi t^{13/6}\sd^3(\tit) 
      \end{bmatrix}$}
    \\\hline\hline
    $\tau_4^a$: & \treeE 
    & $\begin{bmatrix} 0 \\ -6y_2^3y_3-18y_1y_2^2 \end{bmatrix}$
    & $\begin{bmatrix} 0 \\ {\cal O}(t^{3/2}) \end{bmatrix}$ 
    & \jnstrut{$\begin{bmatrix} 
        -\frac{4096}{14553}\sqrt{2}\,c_1^6t^{7/2}\sd'(\tit) \\[\jot] 
        \frac{8192}{43659}\sqrt{2}\,c_1^7t^{23/6}\sd^3(\tit) 
      \end{bmatrix}$}
    \\\hline
    $\tau_4^b$: & \treeF 
    & $\begin{bmatrix} 0 \\ 6y_1^4y_2y_3^2+3y_1^5y_3 \end{bmatrix}$
    & $\begin{bmatrix} 0 \\ {\cal O}(t^{3/2}) \end{bmatrix}$ 
    & \jnstrut{$\begin{bmatrix} 
        \frac{4096}{43659}\sqrt{2}\,c_1^6t^{7/2}\sd'(\tit) \\[\jot] 
        -\frac{8192}{130977}\sqrt{2}\,c_1^7t^{23/6}\sd^3(\tit) 
      \end{bmatrix}$}
    \\\hline
    $\tau_4^c$: & \treeH
    & $\begin{bmatrix} -6y_1y_2^2y_3-6y_1^2y_2 \\ 0 \end{bmatrix}$
    & $\begin{bmatrix} {\cal O}(t^{7/6}) \\ 0 \end{bmatrix}$ 
    & \jnstrut{$\begin{bmatrix} 
        -\frac{4096}{43659}\sqrt{2}\,c_1^6t^{7/2}\sd'(\tit) \\[\jot] 
        \frac{8192}{130977}\sqrt{2}\,c_1^7t^{23/6}\sd^3(\tit) 
      \end{bmatrix}$}
    \\\hline
    $\blt_4$: & \treeG 
    & $\begin{bmatrix} 3y_1^5y_3^2 \\ 9y_1^4y_2y_3^2+3y_1^5y_3 \end{bmatrix}$ 
    & $\begin{bmatrix} {\cal O}(t^{7/6}) \\ {\cal O}(t^{3/2}) \end{bmatrix}$
    & \jnstrut{$\begin{bmatrix} 
        \frac{16384}{43659}\sqrt{2}\,c_1^6t^{7/2}\sd'(\tit) \\[\jot] 
        -\frac{32768}{130977}\sqrt{2}\,c_1^7t^{23/6}\sd^3(\tit) 
      \end{bmatrix}$}
    \\\hline
  \end{tabular}\end{center}
  \caption{Trees of order${}\le4$, with their elementary differentials and integrals for the Emden--Fowler
    equation~\eqref{e:cairy}. The third component is suppressed as it is always zero. In the last column, $\chi = \frac1{4K}
    \int_0^{4K} \sd^2(s) \,\d{s}$ where $4K$ is the period of the function~$\sd$, $\tit=c_1t^{4/3}+c_2$, and only the term of
    leading order is displayed.} 
  \label{t:trees}
\end{table}

\begin{figure}
  \begin{center}
    \includegraphics[width=\linewidth]{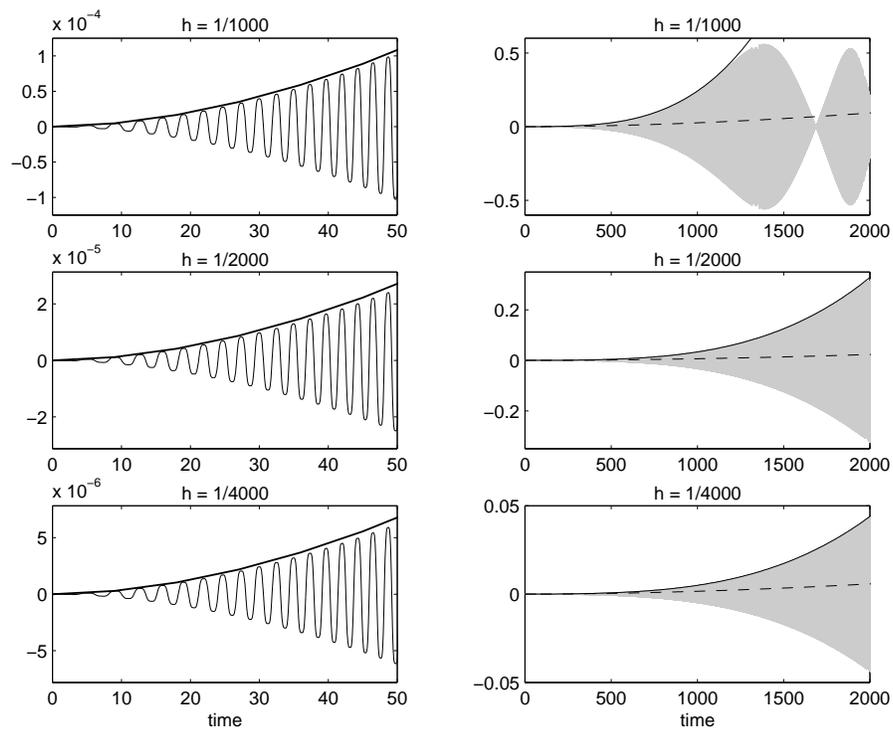}
  \end{center}
  \caption{The first component of the global error committed by Runge's second-order method, and the estimate~\eqref{e:r2-est}.}
  \label{fig3}
\end{figure}

Our first example concerns Runge's second order method, given by
\begin{equation}
  \label{e:runge2}
  \by_{n+1} = \by_n + h\bf \bigl( t_n + \tfrac12h, \by_n + \tfrac12 h\bf(t_n,\by_n) \bigr).
\end{equation}
This is a Runge--Kutta method, hence it can be written as a B-series $B(a_R,\by)$. Some values of the function~$a_R$ are
\begin{gather*}
  a_R\ptreeA=1 ,\quad a_R\ptreeB=1 ,\quad a_R\ptreeC = \tfrac34 ,\quad a_R\ptreeE = \tfrac12, \\
  a_R\ptreeD = a_R\ptreeF = a_R\ptreeG = a_R\ptreeH = 0.
\end{gather*}
It follows from~\eqref{e:b-def} that the B-series coefficients of the modified equation are
\begin{equation*}\label{e:modeqn_airy2}
  \begin{aligned}
    {}&b_R\ptreeA = 1,\quad && b_R\ptreeB = 0,\quad && b_R\ptreeC = -\tfrac14,\quad && b_R\ptreeD = -1, \\
    {}&b_R\ptreeE = 0,\quad && b_R\ptreeF = 0,\quad && b_R\ptreeG = 3,\quad         && b_R\ptreeH = \tfrac32.
  \end{aligned}
\end{equation*}
If we use these numbers in the error estimate~\eqref{e:prop3}, together with the elementary integrals from Table~\ref{t:trees},
we obtain the following estimate for the global error committed by Runge's method:
\begin{equation}
  \label{e:r2-est}
  \bE_h(t) \approx h^2 \begin{bmatrix}
    \frac4{15}\sqrt2\,c_1^4\chi t^{11/6}\sd'(\tit) \\[\jot]
    -\frac{8}{45}\sqrt2\,c_1^5\chi t^{13/6}\sd^3(\tit) 
  \end{bmatrix} 
  +  h^3 \begin{bmatrix}
    \frac{256}{6237}\sqrt2\,c_1^6t^{7/2}\sd'(\tit) \\[\jot]
    -\frac{512}{18711}\sqrt2\,c_1^7t^{23/6}\sd^3(\tit)
  \end{bmatrix}.
\end{equation}
To check this estimate, we compute the solution of the nonlinear oscillator~\eqref{e:cairy} with Runge's second-order
method~\eqref{e:runge2}. The initial conditions are $y(0)=1$ and $y'(0)=0$, which lead to a solution with $c_1\approx0.7$. The
solution is compared to the result of the standard fourth-order Runge--Kutta method with step size $h=1/10000$. According to
Theorem~\ref{th:ef}, this would give an error of about~$10^{-9}$, so we can consider this to be the exact solution. Hence the
global error~$E_h(t)$ is computed by subtracting the result of Runge's method from the ``exact'' solution. The first component is
depicted in Figure~\ref{fig3}. The left column shows the time interval $[0,50]$, and on the right the larger interval $[0,2000]$
is displayed.

In the left column of Figure~\ref{fig3}, the oscillating curve shows the global error~$E_h(t)$. The oscillating behaviour is
correctly predicted by the error estimate~\eqref{e:r2-est}; this is the factor~$\sd'(\tit)$. The amplitude of the oscillations, as
predicted by~\eqref{e:r2-est}, is shown by the thick curve in the left-hand column in Figure~\ref{fig3}. We see that the
estimate~\eqref{e:r2-est} describes the actual error accurately.

The right-hand column of Figure~\ref{fig3} shows a much larger time interval. Here the oscillations are compressed so heavily that
the error appears as a grey blob.  The dashed curve shows the first, leading term of the error estimate~\eqref{e:r2-est}, and the
solid curve shows the sum of both terms. We conclude that the leading $h^2$~term of the estimate does not describe the actual
error correctly, but that the error is predicted accurately if the $h^3$~term is included. For $h=1/1000$ the latter estimate
breaks down around $t=1200$. We note that at that point, the amplitude of the solution is approximately~$0.3$ and the error has
about the same size, so the numerical solution has deviated considerably from the exact solution. For smaller values of the step
size, the estimate~\eqref{e:r2-est} is accurate over the entire time interval $[0,2000]$.

\begin{figure}
  \begin{center}
    \includegraphics[width=\linewidth]{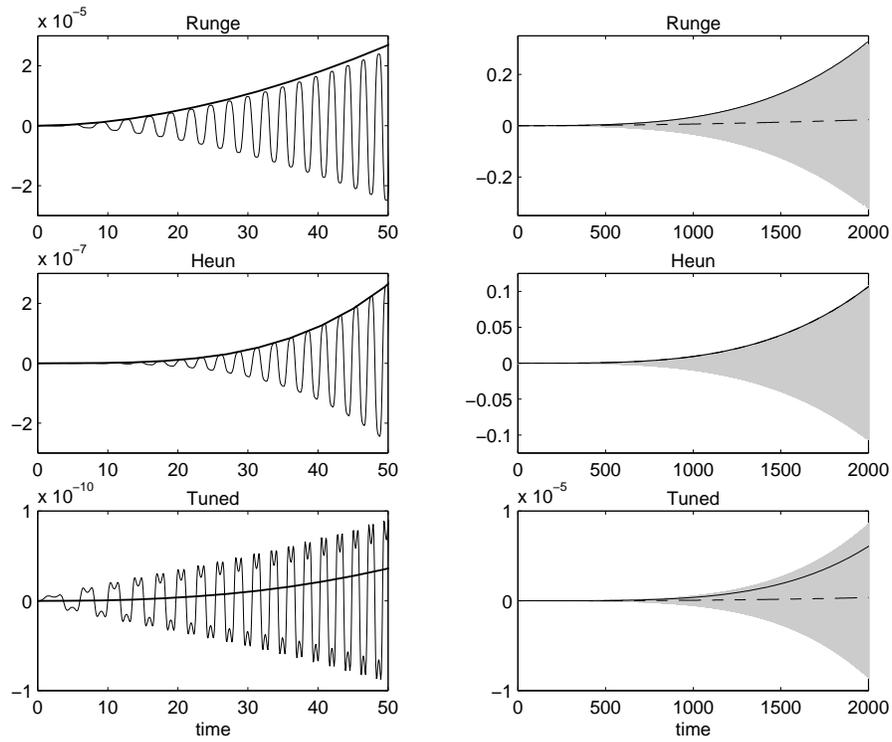}
  \end{center}
  \caption{The first component of the global error committed by Runge's second-order method~\eqref{e:runge2}, Heun's third-order
    method~\eqref{e:heun}, and the specially tuned third-order method~\eqref{e:tuned}, all with step size $h=1/2000$, together
    with their respective error estimates~\eqref{e:r2-est}, \eqref{e:heun-est}, and~\eqref{e:tuned-est}. The lines in the plot
    have the same meaning as in Figure~\ref{fig3}.}
  \label{fig4}
\end{figure}

It follows that for large~$t$ (which here means: in the order of~$1000$), Runge's method behaves essentially as a third-order
method. To check this, we compare it with Heun's classical third-order method:
\begin{equation}
  \label{e:heun}
  \begin{aligned} 
    \bk_1 &= \bf(t_n, \by_n) \\ 
    \bk_2 &= \bf(t_n + \tfrac13h, \by_n + \tfrac13h\bk_1) \\ 
    \bk_3 &= \bf(t_n + \tfrac23h, \by_n + \tfrac13h\bk_2) \\
    \by_{n+1} &= \by_n + h(\tfrac13\bk_1 + \tfrac23\bk_3)
  \end{aligned}
  \qquad\qquad
  \begin{gathered}
    \text{Butcher tableau:} \\
    \begin{array}{c|ccc} 0 \\ 1/3 & 1/3 \\ 2/3 & 0 & 2/3 \\ \hline & 1/4 & 0 & 3/4 \end{array}
  \end{gathered}
\end{equation}
A similar calculation as for Runge's second-order method gives the following estimate for the global error committed by this
method:
\begin{multline}
  \label{e:heun-est}
  \bE_h(t) \approx h^3 \begin{bmatrix}
    -\frac{512}{35721}\sqrt2\,c_1^6t^{7/2}\sd'(\tit) \\[\jot]
    \frac{1024}{107163}\sqrt2\,c_1^7t^{23/6}\sd^3(\tit) 
  \end{bmatrix}  
  +  h^4 \begin{bmatrix}
    \frac{50208}{229635}\sqrt2\,c_1^6t^{5/2}\sd'(\tit) \\[\jot]
    -\frac{60416}{688905}\sqrt2\,c_1^7t^{17/6}\sd^3(\tit)
  \end{bmatrix} \\
  +  h^5 \begin{bmatrix}
    \frac{557056}{34543665}\sqrt2\,c_1^8\chi t^{25/6}\sd'(\tit) \\[\jot]
    -\frac{1114112}{103630995}\sqrt2\,c_1^9\chi t^{9/2}\sd^3(\tit)
  \end{bmatrix}.
\end{multline}
The actual error and the above estimates, for step size $h=1/2000$, are displayed in the second row of Figure~\ref{fig4}. The
first row displays Runge's second-order method. We see that the error estimate~\eqref{e:heun-est} again provides an excellent
description of the actual error. Furthermore, the difference in order between Runge's and Heun's methods shows clearly for small
values of~$t$ (see the left-hand column in Figure~\ref{fig4}). For large values of~$t$ however (cf.\ the right-hand column),
Runge's method behaves essentially as a third-order method and we see indeed that the difference between the two methods is much
smaller.

The bottom row in Figure~\ref{fig4} shows a specially tuned method. Remember from Section~\ref{s:airy} that for the Airy
equation~\eqref{e:airy}, the contribution of the tree~$\blt_k$ to the global error dominates, and this is used by
Orel~\cite{orel:runge-kutta} to construct specially tuned Runge--Kutta methods. We found in Section~\ref{s:global} that this is
not the case for the Emden--Fowler oscillator~\eqref{e:ef}. However, if we carefully study the elementary integrals~${\cal
  I}_\tau$ in Table~\ref{t:trees}, we see that they are scalar multiples of each other:
$$
4 {\cal I}_{\tau_4^a} (t) = -12 {\cal I}_{\tau_4^b} (t) = 12 {\cal I}_{\tau_4^c} (t) = -3 {\cal I}_{\blt_4} (t).
$$
Therefore, bearing in mind the factor~$\alpha(\tau)$ in Theorem~\ref{prop3}, the $h^3$~term will be killed for a third-order
method with
\begin{equation*}
  \label{e:order4cn}
  3b\ptreeE - 3b\ptreeF + b\ptreeH - 4b\ptreeG = 0.
\end{equation*}
From the general theory, it follows that a method has order three if
\begin{equation*}
  \label{e:order3cn}
  b\ptreeA = 1 \quad\text{and}\quad b\ptreeB = b\ptreeC = b\ptreeD = 0.
\end{equation*}
These are five conditions together. An explicit 3-stage Runge--Kutta method has six free parameters, so there is some hope that we
can find such a method satisfying these five conditions. Indeed, it turns out that there exists a one-parameter family of these
methods.  A particular instance is:
\begin{equation}
  \label{e:tuned}
  \begin{aligned} 
    \bk_1 &= \bf(t_n, \by_n) \\ 
    \bk_2 &= \bf(t_n + h, \by_n + h\bk_1) \\ 
    \bk_3 &= \bf(t_n + \tfrac32h, \by_n + \tfrac94h\bk_1 - \tfrac34h\bk_2) \\
    \by_{n+1} &= \by_n + h(\tfrac7{18}\bk_1 + \tfrac56\bk_2 - \tfrac29\bk_3)
  \end{aligned}
  \qquad
  \begin{gathered}
    \text{Butcher tableau:} \\
    \begin{array}{c|ccc} 0 \\ 1 & 1 \\ 3/2 & 9/4 & -3/4 \\ \hline & 7/18 & 5/6 & -2/9 \end{array}
  \end{gathered}
\end{equation}
If we compute the estimate for the global error of this method, we find
\begin{equation}
  \label{e:tuned-est}
  \bE_h(t) \approx h^4 \begin{bmatrix}
    -\frac{5008}{25515}\sqrt2\,c_1^6t^{5/2}\sd'(\tit) \\[\jot]
    \frac{10016}{76545}\sqrt2\,c_1^7t^{17/6}\sd^3(\tit)
  \end{bmatrix} 
  +  h^5 \begin{bmatrix}
    -\frac{78848}{1279395}\sqrt2\,c_1^8\chi t^{25/6}\sd'(\tit) \\[\jot]
    \frac{157696}{3838185}\sqrt2\,c_1^9\chi t^{9/2}\sd^3(\tit)
  \end{bmatrix}.
\end{equation}
As expected, the $h^3$~term has disappeared. This estimate, together with the actual error, is displayed in the bottom row of
Figure~\ref{fig4}. We see that the global error of this method is much smaller than that of Heun's method, though both methods are
of third order. The other thing to note is that the error estimate~\eqref{e:tuned-est} is not as accurate as for the other
methods, even though it predicts the right order of magnitude.

We stress that we are not recommending the use of the method~\eqref{e:tuned}; it works well for \emph{this} equation but we know
of no reason whatsoever why it should also perform well for other equations. It may however be important to gain a better
understanding of the mechanisms that cause this method to work so well.

The reader should keep in mind that the \emph{local} error of the method~\eqref{e:tuned} is still~${\cal O}(h^4)$, as for all
third-order methods. In fact, the B-series coefficient function~$a_T$ of the method~\eqref{e:tuned} satisfies
$$
a_T\ptreeE=\tfrac13 ,\, a_T\ptreeF=2 ,\, a_T\ptreeH=2 ,\, a_T\ptreeG=0.
$$
So the local error of this method is
$$
B(a_T-e, \by) = h^4 \begin{bmatrix} 
  -\frac18y_1^5y_3^2 - \frac14y_1y_2^2y_3 - \frac14y_1^2y_2 \\
  \frac38y_1^4y_2y_3^2 + \frac14y_1^5y_3 + \frac16y_2^3y_3 + \frac12y_1y_2^2
\end{bmatrix} + {\cal O}(h^5). 
$$
where $e$ denotes the function with $e(\emptyset)=0$ and $e(\tau)=1$ for $\tau\neq\emptyset$ (the B-series $B(e,\by)$ gives the
exact solution~\cite{hairer:numerical}, so $B(a_T-e,\by)$ is the local error).

It follows that the \emph{global} error~$\bE_h(t)$ must be~${\cal O}(h^3)$, in contrast to the error estimate~\eqref{e:tuned-est}.
It is only the term of order $h^3 [t^{7/2}, t^{23/6}]^\top$, which for generic third-order methods dominates the $h^3$~term, that
disappears. But there still is an $h^3$~contribution to the global error. We believe that this contribution explains the
substantial difference between the estimated and the actual error for this method, as displayed in Figure~\ref{fig4}. The
difference in the shape of the oscillations of the global error for the last method supports this conclusion.

\section{Discussion}
\label{s:disc}

We used a combination of the Alekseev--Gr{\"o}bner lemma, the theory of modified equations, and asymptotics to calculate an
estimate for the global error committed by a Runge--Kutta method for a class of linear and one of nonlinear oscillators. Numerical
experiments show that these estimates are generally accurate.

While in the linear case, the global error behaves more or less as we expected after the analysis by Iserles~\cite{iserles:on*1},
the nonlinear case brings some surprises which indicate that the entire calculation of Section~\ref{s:global} is necessary. The
specially tuned method~\eqref{e:tuned} shows that it is not enough to consider only the local error. It neither suffices to study
just the first term of order~$h^p$; Figure~\ref{fig4} shows that the next term sometimes dominates. Finally, a comparison of
Theorems~\ref{th:linosc} and~\ref{th:ef} shows that linear and nonlinear oscillators behave differently.

Let us highlight these differences. In the linear case, the contribution of one single tree (namely the tall, branchless tree)
dominates the local error. This is no longer the case for the nonlinear oscillator~\eqref{e:ef}. On the other hand, the
contributions of all trees of the same order to the global error are constant multiples of one another (this has been checked for
order${}\le6$). This allows us to kill the leading term of the global error. Experiments show that this strategy indeed works, but
we still do not completely understand what is going on.

The other difference concerns the integral in the error estimate. The integration may bring down the time exponent due to
oscillations in the integrand. For our equation this happens for even powers of~$h$. As a consequence, the lowest order term in
the expansion of the global error in powers of~$h$ does not always reveal the true behaviour of the global error, as the next term
may dominate for sufficiently large values of~$t$.

We stress that we have only considered one particular class of nonlinear oscillators, namely those satisfying the Emden--Fowler
equation~\eqref{e:ef}.  Hence we do not know how general these phenomena are. We plan to try and repeat the above analysis for
other equations. An obvious target is the generalized Emden--Fowler equation, given by $y'' + a(t) y^n=0$ with
$a(t)>0$~\cite{wong:on}. A more ambitious plan is to try and connect our analysis to the recent interest in highly-oscillatory
systems, triggered by Garc\'\i a-Archilla, Sanz-Serna and Skeel~\cite{garca-archilla.sanz-serna.ea:long-time-step}, and Hochbruck
and Lubich~\cite{hochbruck.lubich:gautschi-type}.

A delicate issue is the validity region of the estimates derived in this paper. If $t$ is too large the ${\cal O}(h^{2p})$~term
will probably dominate, rendering the estimates worthless. On the other hand, the estimates of Sections~\ref{s:airy},
\ref{s:global}, and~\ref{s:num} use the asymptotic solution as $t\to\infty$, so we need $t$ to be sufficiently large. In any case,
the step size~$h$ needs to be sufficiently small; but if it is very small only the $h^p$~term will contribute. The final thing to
keep in mind is that the expansion of the modified equation in powers of~$h$ usually diverges. The only definitive statement we
can make is that more research needs to be done on this issue.

This paper has concentrated on Runge--Kutta methods. Of course, there are many other methods. If one knows that the solution will
be oscillatory, one is probably better off with a method that is specifically designed for this case. The reader is referred to
the review paper by Petzold, Jay and Yen~\cite{petzold.jay.ea:numerical} for examples. We hope that the framework of
Section~\ref{s:theory} is sufficiently flexible to include a wide variety of methods. This requires a theory of modified equations
for these methods, as Hairer~\cite{hairer:backward*1} and Faltinsen~\cite{faltinsen:backward} have done for multistep and
Lie-group methods respectively. Another challenge is presented by variable step size methods.

Let us end by comparing the general estimates from Section~\ref{s:theory}, especially Theorem~\ref{prop2}, with some results
from the literature. The estimate derived by Iserles~\cite{iserles:on*1} can be written as
\begin{equation}
  \label{e:iserles}
  \bE_h(t) = \frac1h \int_{t_0}^t D\Phi_s^t(\by(s)) \, \bL_h\bigl(s,\by(s)\bigr) \, \d{s} + {\cal O}(h^{p+1}),
\end{equation}
where $\bL_h(t,\by)$ denotes the local error of the method. Indeed, this estimate follows from~\eqref{e:prop2}. We note that
\eqref{e:iserles} gives only the ${\cal O}(h^p)$~term of the global error, and we know that this does not always paint an
accurate picture. On the other hand, the local error is easier to determine than the modified equation.

Earlier, Viswanath~\cite{viswanath:global} took quite a different approach but derived a bound on the global error which can be
brought to almost the same form as~\eqref{e:iserles}, namely
$$
\bE_h(t) \le \int_{t_0}^t \| D\Phi_s^t(\by(s)) \| \,\d{s} \cdot \frac1h \max_t |\bL_h(t,\by(t))| + {\cal O}(h^{p+1}).
$$
This bound severely overestimates the actual error in some circumstances. On the other hand, Viswanath~\cite{viswanath:global}
showed that the first factor can be estimated for some large classes of problems, yielding an upper bound for the global error of
any reasonable method.

The final note should be on the older but still highly relevant work of Hairer and Lubich~\cite{hairer.lubich:asymptotic*1}, who
built upon Gragg's asymptotic expansion~\cite{gragg:repeated}. They expanded the global error in powers of the step size~$h$, and
showed how the terms can be found recursively by solving a differential equation involving the previous term. In principle, all
terms of the asymptotic expansion of the global error can be found in this way. Unfortunately, we are not able to solve the
differential equations that appear when we apply their approach to the nonlinear oscillator~\eqref{e:ef}. The first $p$ terms of
Hairer and Lubich's expansion must obviously equal the integral in~\eqref{e:prop2}. However, the connection has not been found
yet.

\subsection*{Acknowledgements}

This work has greatly benefited from discussions of its author with Chris Budd, Stig Faltinsen, Arieh Iserles, Per Christian Moan,
Matthew Piggott, and Divakar Viswanath. The financial support from the EPSRC, Nuffic, VSB Funds, and others is gratefully
acknowledged.

\bibliographystyle{abbrv}
\bibliography{cairy}

\begin{thebibliography}{10}

\bibitem{abraham.marsden:foundations}
R.~Abraham and J.~Marsden.
\newblock {\em Foundations of Mechanics}.
\newblock Benjamin/Cummings, Reading, MA, second edition, 1978.

\bibitem{benettin.giorgilli:on}
G.~Benettin and A.~Giorgilli.
\newblock On the {H}amiltonian interpolation of near-to-the identity symplectic
  mappings with application to symplectic integration algorithms.
\newblock {\em J.~Stat.\ Phys.}, 74(5/6):1117--1143, 1994.

\bibitem{calvo.hairer:accurate}
M.~P. Calvo and E.~Hairer.
\newblock Accurate long-term integration of dynamical systems.
\newblock {\em Appl.\ Numer.\ Math.}, 18:95--105, 1995.

\bibitem{cano.sanz-serna:error*1}
B.~Cano and J.~M. Sanz-Serna.
\newblock Error growth in the numerical integration of periodic orbits, with
  application to {H}amiltonian and reversible systems.
\newblock {\em SIAM J. Numer.\ Anal.}, 34(4):1391--1417, 1997.

\bibitem{cano.sanz-serna:error}
B.~Cano and J.~M. Sanz-Serna.
\newblock Error growth in the numerical integration of periodic orbits by
  multistep methods, with application to reversible systems.
\newblock {\em IMA J. Numer.\ Anal.}, 18(1):57--75, 1998.

\bibitem{chandrasekhar:introduction}
S.~Chandrasekhar.
\newblock {\em An Introduction to the Study of Stellar Structure}.
\newblock University of Chicago Press, 1939.

\bibitem{faltinsen:backward}
S.~Faltinsen.
\newblock Backward error analysis for {L}ie-group methods.
\newblock {\em BIT}, 40(4):652--670, 2000.

\bibitem{garca-archilla.sanz-serna.ea:long-time-step}
B.~Garc{\'\i}a-Archilla, J.~M. Sanz-Serna, and R.~D. Skeel.
\newblock Long-time-step methods for oscillatory differential equations.
\newblock {\em SIAM J. Sci.\ Comput.}, 20(3):930--963, 1999.

\bibitem{gragg:repeated}
W.~Gragg.
\newblock {\em Repeated Extrapolation to the Limit in the Numerical Solution of
  Ordinary Differential Equations}.
\newblock PhD thesis, University of California, Los Angeles, 1964.

\bibitem{hairer:backward*1}
E.~Hairer.
\newblock Backward error analysis for multistep methods.
\newblock {\em Numer.\ Math.}, 84(2):199--232, 1999.

\bibitem{hairer:numerical}
E.~Hairer.
\newblock Numerical geometric integration.
\newblock Lecture notes, 1999.
\newblock Available at http://www.unige.ch/math/folks/hairer/polycop.html.

\bibitem{hairer.lubich:asymptotic*1}
E.~Hairer and C.~Lubich.
\newblock Asymptotic expansions of the global error of fixed-stepsize methods.
\newblock {\em Numer.\ Math.}, 45:345--360, 1984.

\bibitem{hairer.lubich:life-span}
E.~Hairer and C.~Lubich.
\newblock The life-span of backward error analysis for numerical integrators.
\newblock {\em Numer.\ Math.}, 76:441--462, 1997.

\bibitem{hairer.lubich:asymptotic}
E.~Hairer and C.~Lubich.
\newblock Asymptotic expansions and backward analysis for numerical
  integrators.
\newblock In R.~de~la Llave, L.~Petzold, and J.~Lorenz, editors, {\em Dynamics
  of algorithms (Minneapolis, MN, 1997)}, volume 118 of {\em IMA Vol.\ Math.\
  Appl.}, pages 91--106. Springer, New York, 2000.

\bibitem{hairer.nrsett.ea:solving}
E.~Hairer, S.~N{\o}rsett, and G.~Wanner.
\newblock {\em Solving Ordinary Differential Equations {I}. Nonstiff Problems}.
\newblock Springer-Verlag, Berlin, second edition, 1993.

\bibitem{henrici:discrete}
P.~Henrici.
\newblock {\em Discrete Variable Methods in Ordinary Differential Equations}.
\newblock John Wiley \& Sons, New York, 1962.

\bibitem{hochbruck.lubich:gautschi-type}
M.~Hochbruck and C.~Lubich.
\newblock A {G}autschi-type method for oscillatory second-order differential
  equations.
\newblock {\em Numer.\ Math.}, 83(3):403--426, 1999.

\bibitem{iserles:on*1}
A.~Iserles.
\newblock On the global error of discretization methods for highly-oscillatory
  ordinary differential equations.
\newblock To appear in \emph{BIT}, 2001.

\bibitem{iserles.soderlind:global}
A.~Iserles and G.~S{\"o}derlind.
\newblock Global bounds on numerical error for ordinary differential equations.
\newblock {\em J.~Complexity}, 9:97--112, 1993.

\bibitem{lambert:numerical}
J.~Lambert.
\newblock {\em Numerical Methods for Ordinary Differential Equations: {T}he
  Initial Value Problem}.
\newblock John Wiley \& Sons, Chichester, 1991.

\bibitem{neville:jacobian}
E.~Neville.
\newblock {\em Jacobian Elliptic Functions}.
\newblock Clarendon Press, Oxford, 1944.

\bibitem{olver:asymptotics}
F.~Olver.
\newblock {\em Asymptotics and Special Functions}.
\newblock Academic Press, New York-London, 1974.

\bibitem{orel:runge-kutta}
B.~Orel.
\newblock Runge--{K}utta and {M}agnus methods for oscillatory {ODE}s.
\newblock Talk delivered at SciCADE~'01, Vancouver, 2001.

\bibitem{petzold.jay.ea:numerical}
L.~Petzold, L.~Jay, and J.~Yen.
\newblock Numerical solution of highly oscillatory ordinary differential
  equations.
\newblock {\em Acta Numerica}, 6:437--483, 1997.

\bibitem{skeel:thirteen}
R.~Skeel.
\newblock Thirteen ways to estimate global error.
\newblock {\em Numer.\ Math.}, 48(1):1--20, 1986.

\bibitem{viswanath:global}
D.~Viswanath.
\newblock Global errors of numerical {ODE} solvers and {L}yapunov's theory of
  stability.
\newblock {\em IMA J. Numer.\ Anal.}, 21(1):387--406, 2001.

\bibitem{wong:on}
J.~Wong.
\newblock On the generalized {E}mden--{F}owler equation.
\newblock {\em SIAM Review}, 17(2):339--360, 1975.

\end{thebibliography}

\end{document}